\documentclass[letterpaper, 10 pt, conference]{ieeeconf}  

\IEEEoverridecommandlockouts                              
\usepackage{listings}
\usepackage{xcolor}
\usepackage{array} 
\usepackage{amsmath}
\usepackage{float}
\usepackage{amssymb}
\usepackage{amsthm}
\usepackage{url}
\usepackage{graphicx}
\usepackage{subfigure}
\usepackage{xcolor}
\usepackage{booktabs}
\usepackage{multirow}
\usepackage[T1]{fontenc}
\usepackage{lmodern}
\usepackage{times}  

\newcommand{\R}{\mathbb{R}}
\newcommand{\N}{\mathbb{N}}

\usepackage[colorlinks, allcolors=blue]{hyperref}

\usepackage{algorithm}
\usepackage{algpseudocode}
\newtheorem{lemma}{Lemma}
\newtheorem{example}{Example}
\newtheorem{definition}{Definition}
\newtheorem{theorem}{Theorem}
\newtheorem{proposition}{Proposition}
\newtheorem{remark}{Remark}
\newtheorem{assum}{Assumption}

\title{\LARGE \bf Properties of Fixed Points of Generalised Extra Gradient Methods Applied to Min-Max Problems
}

\author{Amir Ali Farzin, Yuen-Man Pun, Philipp Braun, Iman Shames
\thanks{The authors are with CIICADA LAB, School of Engineering,
        Australian National University, Canberra, ACT 2600, Australia
        {\tt\footnotesize \{amirali.farzin,yuenman.pun, philipp.braun,iman.shames\}@anu.edu.au}}
}

\begin{document}
\maketitle
\textit{}
\thispagestyle{empty}
\pagestyle{empty}
\begin{abstract}
This paper studies 
properties of 
fixed points of
generalised Extra-gradient (GEG) algorithms 
applied to
min-max problems. 
We discuss connections
between saddle points of the objective function of the min-max problem
and GEG fixed points.
We show that, under appropriate step-size selections, the set of 
saddle
points (Nash equilibria) is a subset of stable fixed points of GEG. 
Convergence properties of the GEG algorithm are obtained through a stability analysis of a discrete-time dynamical system. 
The results and benefits when compared to existing methods are illustrated through 
numerical examples.
\end{abstract}

\section{Introduction}
Min-max optimisation problems capture the interplay of two decision makers where one seeks to maximise 
an objective 
function while the other aims to minimise it. Traditionally, these problems have been used to
study decision problems in economics \cite{myerson2013game}, and more recently
in adversarial training \cite{madry2018towards} and multi-agent reinforcement
learning \cite{omidshafiei2017deep}. A typical min-max problem takes the form
\begin{equation}\label{eq:eq20}
	\underset{x \in \mathbb{R}^n}{\min}\;\underset{y \in \mathbb{R}^m}{\max}\;f(x,y),
\end{equation}
where  $f: \mathbb{R}^n\times\mathbb{R}^m\rightarrow\mathbb{R}$ is the objective function.

Iterative gradient-based algorithms such as Gradient Descent Ascent (GDA) methods
\cite{korpelevich1977extragradient,chen1997convergence} and their
successors, 
including
extragradient algorithms (EG)~\cite{korpelevich1977extragradient}, optimistic GDA
(OGDA)~\cite{wei2021last}, and stochastic GDA
(SGDA)~\cite{beznosikov2023stochastic} are the go-to methods of choice
for solving these problems. 
In \cite{daskalakis2018limit}, the properties of fixed-points of GDA and OGDA methods are studied, and it is shown that the saddle points of \eqref{eq:eq20} are a subset of the limit-points of the GDA method, and those are in turn a subset of the limit points of the OGDA method.
However, it is known that the GDA method may converge to limit cycles or may even be divergent \cite{daskalakis2018limit}.
This should not come as a surprise. 
Similar statements are also valid when 
continuous-time saddle point-dynamics 
are used to solve min-max problems \eqref{eq:eq20}, which are only guaranteed to converge under specific additional assumptions (see \cite[Ch.~11.5, p.~176]{goebel2024set}, for example).

This absence of convergence guarantees to a fixed point 
of
the GDA method is the main
motivation
for the interest in the study and development of the Extra Gradient method  and its variants as they are shown to have better convergent properties \cite{korpelevich1977extragradient,chae2023two,diakonikolas2021efficient,kim2024double}. 

Inspired by the results in these papers
and aiming to provide a unified framework for studying EG algorithms, we introduce the \emph{Generalised Extra Gradient (GEG)} algorithm 
\begin{equation}\tag{GEG}
\begin{aligned}\label{dy}
	\hat{x}_k &= x_k -h_{1x} \nabla_x f(x_k,y_k)\\
	\hat{y}_k &= y_k +h_{1y} \nabla_y f(x_k,y_k)\\
	x_{k+1} &= x_k -h_{2x} \nabla_x f(\hat{x}_k,\hat{y}_k)\\
	y_{k+1} &= y_k +h_{2y} \nabla_y f(\hat{x}_k,\hat{y}_k),
\end{aligned}
\end{equation}
where $h_{1x}$, $h_{2x}$, $h_{1y}$, and $h_{2y}$ are positive step-sizes. To simplify presentation, in the rest of this paper, we define $\gamma=\frac{h_{2x}}{h_{1x}}=\frac{h_{2y}}{h_{1y}}$ and $\tau=\frac{h_{1y}}{h_{1x}}=\frac{h_{2y}}{h_{2x}}$.
This algorithm is a generalisation of other  
EG algorithms 
in terms of possible step-size selections.
For example, the two timescale EG algorithm in \cite{chae2023two} is a special case of \eqref{dy} for the case where $\gamma=1$ and $\tau\geq 1$. Different 
variants of  
EG algorithms, characterised by parameters $\tau$ and $\gamma,$ are summarised 
in Table~\ref{tab:egv}.
\begin{table}[thb]
\centering
\begin{tabular}{c|c|c}
\toprule
Algorithm & \multicolumn{2}{c}{Parameters} \\
\midrule
EG~\cite{korpelevich1977extragradient} &  $\tau=1$ & $\gamma=1$  \\
\midrule
$\tau$-EG~\cite{chae2023two}&  $\tau\geq1$ & $\gamma=1$  \\
\midrule
EG+~\cite{diakonikolas2021efficient} &  $\tau=1$ & $0<\gamma\leq 1$  \\
\midrule
\eqref{dy} &  $\tau>0$ & $\gamma>0$  \\
\bottomrule
\end{tabular}
\caption{Variants of extragradient algorithm}
\label{tab:egv}
\end{table}

\vspace{-0.2cm}

In this paper, we use tools from dynamical systems theory to establish 
stability properties of fixed points of \eqref{dy} via studying the equilibria of an equivalent dynamical system. Specifically, leveraging the properties of hyperbolic equilibria, we demonstrate that \eqref{dy} avoids its unstable fixed points for almost all initialisations and show that, with appropriate step-size selection, the set of 
saddle points of \eqref{eq:eq20} is a subset of the stable equilibria of the dynamical system of interest and consequently the stable fixed points of \eqref{dy}.

\textit{Outline:}  The necessary preliminaries as well as the dynamical system representation of \eqref{dy} 
are introduced in Section~\ref{sec:per}. In Section~\ref{sec:main},  the main results on the equilibria of the dynamical system are presented. In Section~\ref{sec:exps}, three illustrative examples are studied. The conclusions and future research directions are discussed in Section~\ref{sec:conc}. 

\textit{Notation:}
For , 
$x,y\in \mathbb{R}^d$, 
$\langle x , y \rangle=x^\top y$ and $\|x\|=\sqrt{\langle x,x\rangle}$. Moreover,  
for $c=a+\jmath b \in\mathbb{C}$,
 $\Re(c)=a$ and $|c|=\sqrt{a^2+b^2}$.
A ball of radius $\delta >0$ around $z^\ast$ is denoted by $\mathcal{B}_\delta(z^\ast)=\{z\in \R^d| \ \|z-z^\ast\|\leq \delta\}$.
For $r,n,m\in \N$, 
$f:\R^{n}\rightarrow \R^m$ is a $C^r$-function if it is $r$-times continuously differentiable.
Gradient and Hessian of 
$f: \mathbb{R}^{n+m} \to \mathbb{R}$ are denoted by 
$\nabla f$ and  $\nabla^2 f$. 
For $(x,y)\in \R^{n+m}$, $ \nabla_x f$, $\nabla_y f$, $\nabla^2_{xx} f$, $ \nabla^2_{xy} f$, and $\nabla^2_{yy} f$ 
denote parts of the gradient/Hessian where the derivative of $f$ is taken with respect to the subscript. 
For $A\in \R^{n\times n}$, the spectrum and spectral radius are denoted by
 $\sigma(A) = \{\lambda \in \mathbb{C} | \ \det(A-\lambda I) = 0 \}$ and $
 \rho(A) = \max_{\lambda \in \sigma(A)} |\lambda|$, respectively. 
For $A$ symmetric, $A \succ 0$ ($A \succeq 0$) indicates that  \( A \) is positive (semi-positive) definite and $I$ denotes the identity matrix.

\section{Preliminaries \& Problem Formulation} \label{sec:per}
In this paper, we investigate relationships between saddle points of the min-max problem \eqref{eq:eq20} and fixed points of \eqref{dy}. Instead of investigating  \eqref{dy} directly, we discuss the algorithm through the lens of stability properties of discrete-time dynamical systems ($w: \R^{d} \rightarrow \R^d$)
\begin{align}
    z_{k+1} = w(z_k) \qquad (\text{or } z^+=w(z))  \label{eq:discrete_time_sys}
\end{align}
with state
$z\in \R^d$. 
\begin{definition}
Consider the discrete-time system \eqref{eq:discrete_time_sys}. 
A point $z^e\in \R^d$ is called 
an equilibrium of $w$ if  $ w(z^e) = z^e$.
\end{definition}
To 
formulate \eqref{dy} as a discrete-time dynamical system, we 
assume that $f$ is $C^3$ and define
\begin{align}
    z=\left[ \hspace{-0.075cm} \begin{array}{c}
         x  \\
         y 
    \end{array} \hspace{-0.075cm} \right],  \ 
    \Lambda_\tau =\left[ \hspace{-0.075cm} \begin{array}{cc}
         \frac{1}{\tau} I & 0  \\
         0 & I 
    \end{array} \hspace{-0.075cm} \right], \
    F(z)=\left[ \hspace{-0.075cm} \begin{array}{r}
         \nabla_x f(z)\\ 
         -\nabla_y f(z)
    \end{array} \hspace{-0.075cm} \right]. \label{eq:F_definition}
\end{align}
Using the definitions~\eqref{eq:F_definition} and $h_{1,y}=\eta$, and by eliminating the variables $\hat{x}_k$, $\hat{y}_k$, the GEG algorithm can be written as 
\begin{align}\label{dy2}
    z_{k+1} = w(z_k) = z_k - \gamma\eta\Lambda_\tau F(z_k-\eta\Lambda_\tau F(z_k)).
\end{align}
It can be observed that  
$(\hat{x},\hat{y},x,y)$ is a fixed point of the GEG algorithm 
if and only if $\nabla f(x,y)=0$. Equivalently, a fixed point of \eqref{dy}
is characterised through  an equilibrium $z^e$ of the  
dynamics \eqref{dy2} where $F(z^e)=0$ is satisfied.
With  
\begin{align}
H(z) &= \begin{bmatrix}
    -\nabla_{xx} f(z)& -\nabla_{xy} f(z)\\ \nabla_{yx} f(z)& \nabla_{yy} f(z)
\end{bmatrix}  =
\left[ \begin{smallmatrix}
    -I & 0  \\
    0 & I
\end{smallmatrix} \right] \nabla^2 f(z)
\label{eq:DF_def}
\end{align} 
denoting the  
derivative of $-F$, 
the Jacobian of 
\eqref{dy2}
at $z$ is
\begin{align}\label{JGEG}
    J(z) = \mathrm{I}  +  \gamma\eta\Lambda_\tau H(z-\eta\Lambda_\tau F(z))\big(\mathrm{I}+\eta\Lambda_\tau H(z)\big).
\end{align}

To analyse stability properties of equilibria of the dynamical system \eqref{eq:discrete_time_sys}, we use standard  
stability definitions and corresponding results in
\cite{hespanha2018linear,kellett2023introduction}.

\begin{definition}[\mbox{\cite[Def. 5.2]{kellett2023introduction}}]
    Consider an equilibrium $z^e\in \R^d$ of 
    \eqref{eq:discrete_time_sys}. Then $z^e$ 
    is  
    stable  
    if, for any $\varepsilon >0$ there exists $\delta >0$ such that if $z_0 \in \mathcal{B}_\delta(z^e)$ then, for all $k\geq 0$, it holds that $ z_k \in \mathcal{B}_\varepsilon(z^e).$
    The 
    equilibrium is asymptotically stable if it is stable and there exists $\delta>0$ such that if $z_0 \in \mathcal{B}_\delta (z^e)$ then $\lim_{k\rightarrow \infty} z_k = z^e$.
    If $z^e$ 
    is not stable, then it is called unstable.
\end{definition}

Note that asymptotic stability of an equilibrium of \eqref{dy2} is equivalent to stability of a critical point of a fixed point iteration such as \eqref{dy}. Similarly, instability with respect to a dynamical system is equivalent to instability of a critical point in a fixed point method.

\begin{theorem}[\mbox{\cite[Thm 8.7]{hespanha2018linear}}] \label{thm:stability_of_equilibria}
   Let $z^e \in \R^n$ denote an equilibrium of 
   \eqref{eq:discrete_time_sys} with $C^2$ function $w(\cdot)$. If  the spectral radius of $J(z^e)$ 
   satisfies $|\rho(J(z^e))| <1$, then $z^e$ is locally asymptotically stable. If $|\rho(J(z^e))| >1$, then $z^e$ is unstable.
\end{theorem}

\begin{remark}
For $|\rho(J(z^e))| =1$, no statement about stability/instability of the equilibrium $z^e$ can be made, despite the contradicting assertion in \cite{daskalakis2018limit}. 
Thus, we deliberately deviate from the presentation of stability notions as given in \cite{daskalakis2018limit}.
\end{remark}

Stability properties of nonlinear systems are in general only local. 
For an equilibrium $z^e\in \R^d$ of the discrete-time system \eqref{eq:discrete_time_sys}, we define its region of attraction as the set
\begin{align}
    \mathcal{R}(z^e) = \{z_0 \in \R^d| \lim_{k\rightarrow \infty} z_k = z^e, \ z_{k+1}=w(z_k)\}. \label{eq:RoA}
\end{align}

As a next step, we recall
standard terminology for 
min-max problems \eqref{eq:eq20}.

\begin{definition}[Critical point {\cite[Def. 1.6]{daskalakis2018limit}}] Let $f: \mathbb{R}^{n+m}
\rightarrow \mathbb{R}$ be a 
$C^1$ function.
 A point $(x^\ast,y^\ast)\in \R^{n+m}$ is called a critical point of $f$ if $\nabla f(x^\ast,y^\ast)=0$.
\end{definition}

\begin{definition}[  
{\cite[Def. 3.4.1]{bertsekas2009convex}}]
    A point $(x^\ast,y^\ast)\in \R^{n+m}$ is a (local) saddle point of $f:\mathbb{R}^{n+m} \rightarrow \mathbb{R}$ 
    if there exists a neighbourhood $U$ around $(x^\ast,y^\ast)$ so that for all $(x,y)\in U$, 
	\begin{align}
    f(x^\ast,y)\leq f(x^\ast,y^\ast)\leq f(x,y^\ast). \label{eq:saddle_local}
	\end{align}
        If inequality 
    \eqref{eq:saddle_local} is 
    strict, the 
    saddle 
    point is 
    strict.
\end{definition}

In the context of min-max problems, saddle points are equivalently referred to as 
Nash equilibria. Throughout this paper, we limit ourselves to using the saddle points terminology.

\begin{proposition}[{\cite[Prop. 4 and 5]{jin2020local}}]\label{locNash}
    Let $f: \mathbb{R}^{n+m}
    \rightarrow \mathbb{R}$ be a $C^2$ function. Then
    any 
    saddle point
    $(x,y)$ is a critical point of $f$ and satisfies $\nabla_{xx}^2 f(x,y) \succeq 0$ and $\nabla_{yy}^2 f(x,y) \preceq 0$. 
    It is a strict  
    saddle point if $\nabla_{xx}^2 f(x,y) \succ 0$ and $\nabla_{yy}^2 f(x,y) \prec 0.$
\end{proposition}

In the following we investigate connections between local asymptotic stability of equilibria of \eqref{dy2} and saddle points of  
$f$. 
To this end, we rely on the following assumptions.
\begin{assum}[Invertibility of $\nabla^2 f$]\label{as:inv}
	Let $f: \mathbb{R}^{n+m}
    \rightarrow \mathbb{R}$ be a $C^3$ function and let  
     $(x^*,y^*)\in \R^{n+m}$ be a saddle point.   
    The Hessian $\nabla^2 f(x^*,y^*)$ is invertible. 
\end{assum}
Invertibility of the Hessian of
$f$ is a common assumption in the 
analysis of iterative optimisation algorithms, e.g., see \cite{daskalakis2018limit,chae2023two}. 
We only require
invertibility of the Hessian at  
saddle points 
and not for all $(x,y)\in\mathbb{R}^{n+m}$. 
At strict
saddle
points the Hessian is always invertible \cite{daskalakis2018limit}. 
\begin{assum}\label{as:c2}
Function $f: \mathbb{R}^{n+m}\rightarrow \mathbb{R}$ is a $C^3$ function and $\nabla f$ is globally Lipschitz with constant $L>0.$ 
\end{assum}

\section{Analysis of generalised Extra-gradient} \label{sec:main}

 In this section, we study the local behaviour of \eqref{dy} when used to
 find 
 saddle points  
 of
 \eqref{eq:eq20}. Instead of studying \eqref{dy} directly, we study 
 the dynamical system~\eqref{dy2}.
 We first present three auxiliary results used to derive the main results. The main results of the paper are presented in Section \ref{sec:main_results}.

 \subsection{Derivation of preparatory and auxiliary results}

 In the following, we show that under appropriate parameter selection
 , $w(\cdot)$ in \eqref{dy2} defines
 a local diffeomorphism\footnote{A local diffeomorphism is a function that is
 locally invertible, smooth, and has a smooth local inverse~\cite{lee2012smooth}.}.

\begin{lemma}\label{lem:defeo}
	Under Assumption~\ref{as:c2} consider $w$  
    in 
    \eqref{dy2} and let  
 $h_{1y} =\eta$ 
 and $\eta\in(0,\frac{c}{L})$ for some positive constant $c>0$. 
 If
 \begin{align}\label{paramchoice}
     (\tau,\gamma) \in \mathcal{P}_1 \cup \mathcal{P}_2,
 \end{align}
 where $\mathcal{P}_1 =  \{ (\tau,\gamma)| 0<\tau\leq1,\ 0< \gamma \leq \frac{\tau^2}{c\tau+c^2} \}$ and 
 $\mathcal{P}_2 = \{ (\tau,\gamma)|
    \tau\geq 1,\ 0< \gamma \leq \frac{1}{c+c^2}\}$,
then $w(\cdot)$ defines a local diffeomorphism.
\end{lemma}

\begin{proof}
	Considering the Inverse Function theorem \cite[Thm. 9.24]{rudin1964principles}, it is sufficient to show that $J$ in \eqref{JGEG}
is invertible.
 To this end, we 
 show that none of the eigenvalues of $J$ is zero.
   Let $A = \eta\Lambda_\tau H(z-\eta\Lambda_\tau F(z))$  and $B= \eta\Lambda_\tau H(z)$:
 \begin{align}
     \|\gamma A (\mathrm{I}+B)\| \leq \gamma \|A \| \|I+B \| \leq \gamma\|A\|(1+\|B\|). \label{eq:estProofLem1}
 \end{align}
Recalling the definition of $\Lambda_\tau$ in \eqref{eq:F_definition}, the definition of $H$ in \eqref{eq:DF_def} and using the fact that $\nabla f$ 
is globally Lipschitz
by assumption, the estimates
\begin{align}
\begin{split}
    \|A\| &\leq \eta \| \Lambda_\tau\| \|H(z-\eta\Lambda_\tau F(z))\| \leq \eta  \max\{\tfrac{1}{\tau},1\} L, \\
    \|B\| &\leq \eta \|\Lambda_\tau\| \|H(z)\| \leq \eta \max\{\tfrac{1}{\tau},1\} L 
\end{split} \label{eq:global_L_step}
\end{align} 
hold. Using these expressions in \eqref{eq:estProofLem1} together with $\eta\leq \frac{c}{L}$:
\begin{align*}
     \|\gamma A (\mathrm{I}+B)\| &\leq \gamma \eta L \max\{\tau^{-1},1\}  (1+ \eta L \max\{\tau^{-1},1\} ) \\
     &\leq \gamma    (c\max\{\tau^{-1},1\}+ c^2  \max\{\tau^{-2},1\} ).
 \end{align*}
The choice of $\gamma$ satisfying $\eqref{paramchoice}$ leads to $\|\gamma A (\mathrm{I}+B)\| <1$. This in turn results in invertibility of $J=I+\eta A(I+B)$.
\end{proof}

\begin{remark} \label{rem:Lipschitz}
    Lemma \ref{lem:defeo} requires the gradient of $f$ to be globally Lipschitz for the estimates in \eqref{eq:global_L_step}. If the domain of interest is restricted to a compact domain, local Lipschitz properties in Lemma \ref{lem:defeo} and in the following are sufficient.
\end{remark}

Next, we show that the eigenvalues of $\Lambda_\tau H$ for $\Lambda_\tau$ and $H$ defined in \eqref{eq:F_definition} and \eqref{eq:DF_def}, respectively, have non-positive real parts  at local saddle points for any value of $\tau>0$. 
\begin{lemma}\label{negeig}
    Let $f: \mathbb{R}^{n+m} 
    \rightarrow \mathbb{R}$ be a $C^2$ function, $H$ and $\Lambda_\tau$ be defined in \eqref{eq:DF_def} and \eqref{eq:F_definition}. Let $z^\ast=(x^\ast,y^\ast)$ be a
    saddle point of $f$ and $\kappa\in\sigma(\Lambda_\tau H(z^\ast))$. 
    Then $\Re(\kappa)\leq0.$
\end{lemma}

\begin{proof}
We use the notation $\omega = 1/\tau$, define the matrix
        $ U = \left[\begin{smallmatrix}
		0 & \sqrt{\omega}\mathrm{I}\\
		\mathrm{I} & 0 
        \end{smallmatrix} \right]$    

        and observe that $\nabla_{xy}^2 f=(\nabla_{yx}^2 f)^\top$ since $f$ is a $C^2$ function.     
	Dropping the arguments of the functions to shorten the expressions in the following, it holds that   
    \begin{align*}
    \bar{H}_\tau = U(\Lambda_\tau H)U^{-1}= \begin{bmatrix}
		\nabla_{yy}^2f & \sqrt{\omega} (\nabla_{xy}^2 f)^\top \\
		-\sqrt{\omega} \nabla_{xy}^2 f & -\omega\nabla_{xx}^2 f
        \end{bmatrix}
        \end{align*} 
and  $\Lambda_\tau H$ and $\bar{H}_\tau$ have the same eigenvalues since $U$ is invertible. 
    Consider 
    $
        \tfrac{1}{2}(\bar{H}_\tau+\bar{H}^\top_\tau) = 
        \left[ \begin{smallmatrix}
		\nabla_{yy}^2f &0 \\
		0 & -\omega\nabla_{xx}^2f 
        \end{smallmatrix} \right].
    $
    From Proposition~\ref{locNash}, at a saddle point, we have \(\nabla_{xx}^2 f \succeq 0\) and \(\nabla_{yy}^2 f \preceq 0\). Consequently, \(\frac{1}{2}(\bar{H}_\tau+\bar{H}_\tau^\top)\) is negative semi-definite.  
Applying Ky Fan’s inequality \cite[Page 1315]{moslehian2012ky}, which states that the sequence of eigenvalues (in decreasing order) of $\frac{1}{2}(\bar{H}_\tau+\bar{H}_\tau^\top)$ majorizes the real part of the eigenvalues of $\bar{H}_\tau$, we conclude that for any \(\lambda\in\sigma(\bar{H}_\tau)\),  
\begin{align*}
\Re(\lambda) \leq \tfrac{1}{2} \lambda_{\max}(\bar{H}_\tau+\bar{H}_\tau^\top) \leq 0.
\end{align*}
This completes the proof as \( \sigma(\Lambda_\tau H) =\sigma (\bar{H}_\tau)\). 
\end{proof}
Next, we show that the eigenvalues of $J,$ 
in \eqref{JGEG}, are related to the eigenvalues of $\Lambda_\tau H$ at a fixed point via a bijection. 
\begin{lemma}\label{Jeig}
	Let $f: \mathbb{R}^{n+m} 
    \rightarrow \mathbb{R}$ be a 
    $C^2$-function, $\Lambda_\tau$, $H$, and $J$  be defined in \eqref{eq:F_definition}, \eqref{eq:DF_def}, and \eqref{JGEG}, respectively, and let $z^\ast=(x^\ast,y^\ast)$ be a fixed point of $f$. Moreover, let $\kappa\in\sigma(\Lambda_\tau H(z^\ast))$, 
    $h_{1y}=\eta,$ $\tau>0$, $\gamma>0$, $\xi = 1+\eta\kappa,$ and $\beta = \gamma\eta\kappa$. Then for any $\lambda\in\sigma(J(z^\ast))$, we have
    \begin{align}\label{eig}
		\lambda =  1\hspace{-0.05cm}+\hspace{-0.05cm}\xi\beta = 1\hspace{-0.05cm}+\hspace{-0.05cm} (1 + \eta \kappa) \gamma \eta \kappa = 1 + \gamma (\eta \kappa + \eta^2 \kappa^2).	\end{align} 
\end{lemma}
\begin{proof}

Since $z^\ast$ is a critical point, it holds that 
$F(z^\ast)=0$ 
and $J(z^\ast)$ can be written as
\begin{align*}
J(z^\ast) = \mathrm{I} +\gamma\eta\Lambda_\tau H(z^\ast)\big(\mathrm{I}+\eta\Lambda_\tau H(z^\ast)\big).
     \end{align*}
    Let $\nu$ denote an eigenvector corresponding to an eigenvalue $\kappa$,  
    i.e., it holds that $\Lambda_\tau H \nu = \kappa \nu$
and  $(\Lambda_\tau H)^2 \nu = \kappa^2 \nu.$ 
In particular, 
the eigenvalues of \((\Lambda_\tau H)^2\) are \(\kappa^2\), with the same eigenvectors as \(\Lambda_\tau H\). 
By direct calculation and by omitting the function arguments, we obtain the following relationship  
\begin{align*}
    (\gamma\eta\Lambda_\tau H &+ \gamma\eta^2(\Lambda_\tau H)^2)\nu 
    = \gamma\eta^2(\Lambda_\tau H)^2\nu+\gamma\eta\Lambda_\tau H\nu \\ 
    &= \gamma\eta^2\kappa^2\nu+\gamma\eta\kappa\nu  
    = ( \gamma\eta\kappa+ \gamma\eta^2\kappa^2)\nu.
\end{align*}  
This implies that \(\gamma\eta\kappa+ \gamma\eta^2\kappa^2\) is an eigenvalue of \(\gamma\eta\Lambda_\tau H + \gamma\eta^2(\Lambda_\tau H)^2\) and hence 
$
\lambda = 1 + \gamma\eta\kappa+ \gamma\eta^2\kappa^2 = 1 + (\eta\kappa+1)\gamma\eta\kappa= 1 + \xi\beta.
$  
\end{proof}

\subsection{Saddle points and GEG stable points} \label{sec:main_results}
In this section, we 
derive properties of saddle points of the min-max problem \eqref{eq:eq20} through stability properties of equilibria of the dynamical system \eqref{dy2}. First, we show that the region of attraction of an unstable equilibrium of \eqref{dy2} is of measure zero.

\begin{theorem}\label{nouns}
	Let $f$ satisfy Assumption~\ref{as:c2}, $w$ be defined according to 
    \eqref{dy2}, and 
    let  $c>0$, 
    $\gamma>0,$ 
    $h_{1y} = \eta$ and $\eta\in(0,\frac{c}{L}).$ If \eqref{paramchoice} holds, then
    for any unstable equilibrium $z^\ast=(x^\ast,y^\ast)$ of \eqref{dy2}, the region of attraction $\mathcal{R}(z^\ast)\subset \R^d$ is of measure zero.

\end{theorem}

\begin{proof}
The proof follows similarly to that of \cite[Thm 2]{lee2019first}. By applying Lemma~\ref{lem:defeo}, we establish that the right-hand side of \eqref{dy2} 
is a diffeomorphism. Then, the remainder of the proof follows along the lines of the arguments in \cite{lee2019first}.
\end{proof}

Theorem \ref{nouns}
shows that if one chooses the initial condition of \eqref{dy2} randomly, the probability of \eqref{dy2} converging to an unstable equilibrium is zero. 

The following results provide 
connections between the set of asymptotically stable equilibria of \eqref{dy2} and the set of saddle points of the min-max problem \eqref{eq:eq20}. 
Theorem \ref{th:real} is restricted to the case that the eigenvalues of $H$ are real while Theorem \ref{th:general} discusses the general case.

\begin{theorem}\label{th:real}
Under Assumptions \ref{as:inv} and \ref{as:c2}, let 
$\tau >0$ and $h_{1y}=\eta\in(0,\frac{\min\{1,\tau\}}{L})$.
For  $\Lambda_\tau$ and $H$  defined in \eqref{eq:F_definition} and \eqref{eq:DF_def}, respectively, assume that $\sigma(\Lambda_\tau H(z^\ast))\subset \R$ for all saddle points 
$z^\ast\in \R^d$ of $f$.
Then, for $\gamma \in (0,8)$, the set of saddle points 
is a subset of the set of asymptotically stable equilibria of \eqref{dy2}. 
\end{theorem}

\begin{proof}
Let  $\zeta=\eta\kappa$ with $\kappa\in \sigma(\Lambda_\tau H(z^\ast))$. 
We note that $\zeta$ is real since $\kappa$ is real by assumption.
Considering Lemma~\ref{Jeig}, the eigenvalues of $J(z^\ast)$ 
are in the form 
$ \lambda = 1+\xi \beta  = 1+\gamma(\zeta+\zeta^2).$
To conclude asymptotic stability from Theorem \ref{thm:stability_of_equilibria}, the condition
$|1+\gamma (\zeta+\zeta^2)|<1$, 
i.e., $\zeta+\zeta^2<0$ and 
$-2<\gamma(\zeta+\zeta^2)$, needs to be satisfied. %
The first condition 
is satisfied for $\zeta\in(-1,0)$ and which implies that
$\zeta+\zeta^2\in 
[-0.25,0)$. Hence, a sufficient condition for  $-2<\gamma(\zeta+\zeta^2)$ to hold for  $\zeta \in(-1,0)$ 
is $\gamma<8$.

From
Lemma~\ref{negeig}  
(and since 
$\kappa$ is real) we know that $\zeta\leq0$ holds. 
According to Assumption~\ref{as:inv}, $\Lambda_\tau H(z^\ast)$ is invertible at 
saddle points 
and thus $\zeta\neq0.$
Using the fact
that $\nabla f$ is $L$-Lipschitz
and 
due to the parameter selection
$\eta$ and $\tau$, it follows 
that  $|\zeta| \leq \rho(\eta\Lambda_\tau H(z^\ast))\leq\|\eta\Lambda_\tau H(z^\ast)\|<1$ as $\|\Lambda_\tau\|
\leq1$ for $\tau\geq1$ and $\|\Lambda_\tau\|\leq\frac{1}{\tau}$ for $0<\tau\leq1.$ Thus, we have shown that $\zeta\in(-1,0)$ and the assertion 
follows.
\end{proof}
Note that if in Theorem~\ref{th:real}, the upper-bound on the step-size $\eta$ is tighter, i.e., $\eta\in(0,c\frac{\min\{1,\tau\}}{L})$, $c\in(0,\frac{1}{2})$, then the upper bound on $\gamma$ can be larger accordingly, i.e., $\gamma \leq \frac{2}{c(1-c)}$. 
\begin{theorem}\label{th:general}
Let Assumptions \ref{as:inv} and \ref{as:c2} be satisfied, let $\tau >0$ and $h_{1y}=\eta\in(0,\frac{\min\{1,\tau\}}{L})$.
    Then, for $\gamma\in (0,1]$, the set of saddle points $z^\ast \in \R^d$
    of $f$ are a subset of the set of asymptotically stable equilibria of \eqref{dy2}.
\end{theorem}

\begin{proof}
    Let
    $\kappa \in \sigma(\Lambda_\tau H(z^\ast))$. Hence, $1+\eta\kappa \in \sigma(I+\eta\Lambda_\tau H(z^\ast))$ and $\gamma\eta\kappa\in\sigma(\gamma\eta\Lambda_\tau H(z^\ast))$. 
From Lemma~\ref{negeig},  we know that $\Re(\eta\kappa)\leq0$. Thus, $a\leq 0$, where $\eta\kappa=a+\jmath b $ for $a,b\in \R$. From
Assumption~\ref{as:inv}, we can  
conclude that $(a,b)\neq (0,0)$ since $\eta\Lambda_\tau H(z^\ast)$ is invertible at saddle points of $f$. 
By the choice of $\eta$ and $\tau$ and the fact that $\|\nabla^2f\|\leq L$, it holds that $\rho(\eta\Lambda_\tau H(z^\ast))  \leq \|\eta\Lambda_\tau H(z^\ast)\| <1$.
Hence, $a\in(-1,0]$ and $b\in(-\sqrt{1-a^2},\sqrt{1-a^2})$.
	From 
    Lemma~\ref{Jeig} eigenvalues of $J(z^\ast)$  
    can be written as 
    \begin{align*}
        \lambda &=1+\gamma\eta\kappa(1+\eta\kappa)
        =1+\gamma(a+\jmath b)(1+a+\jmath b)
        \\&= 1 + \gamma(a^2-b^2+a)+\jmath \gamma(2ab+b).
    \end{align*}
    To establish a
    connection 
    between
    saddle points 
    of $f$ and  
    the set of asymptotically stable equilibria of \eqref{eq:discrete_time_sys}, we derive a
    bound on $\gamma$ to ensure that $|\lambda|<1,$ i.e.,   
    $|\sqrt{(1 + \gamma((a^2-b^2+a))^2+\gamma^2(2ab+b)^2}|<1$ or equivalently $|(1 + \gamma((a^2-b^2+a))^2+\gamma^2(2ab+b)^2|<1$ needs to be satisfied. 
    A
    bound on $\gamma$ ensuring $|\lambda|\leq1$ is obtained through
    Mathematica\footnote{The Mathematica  
    code is reported  
    in Appendix~\ref{mtc}.} and in particular, it holds that
    $|\lambda|\leq1$ if $\gamma\in(0,1]$.  
    To complete the proof,
    we need to analyse the cases 
    where $|\lambda|=1$ is satisfied. 
    The trivial solutions to 
    $(1 + \gamma((a^2-b^2+a))^2+\gamma^2(2ab+b)=1$, or equivalently
    $((a^2 \hspace{-0.05cm} -\hspace{-0.05cm} b^2 \hspace{-0.05cm}+ \hspace{-0.05cm}a)^2 \hspace{-0.05cm}+ \hspace{-0.05cm}(2ab\hspace{-0.05cm}+ \hspace{-0.05cm}b)^2)\gamma^2\hspace{-0.05cm}+\hspace{-0.05cm}(2(a^2\hspace{-0.05cm}-\hspace{-0.05cm}b^2\hspace{-0.05cm}+\hspace{-0.05cm}a))\gamma=0$ 
    are $\gamma=0$ (which is ruled out since $\gamma$ is positive by assumption), $(a,b)=(0,0)$ (which is ruled out by Assumption~\ref{as:inv}), and $(a,b)=(-1,0)$ (which is ruled out since $a\in (-1,0]$).  
    If $\gamma\neq 0$, 
    $\gamma = \frac{-2(a^2-b^2+a)}{(a^2-b^2+a)^2+(2ab+b)^2}.$ Using Mathematica\footnote{The Mathematica
    code is reported
    in Appendix~\ref{mtc2}.}, it can be seen that
    $\frac{-2(a^2-b^2+a)}{(a^2-b^2+a)^2+(2ab+b)^2}>1$ in the feasible range of $a$ and $b$. Choosing a positive $\tau$, $\eta\in(0,\frac{\min\{1,\tau\}}{L})$, and $\gamma \in (0,1],$  then the set of saddle points  
     of $f$ are a subset of the set of asymptotically stable points of \eqref{dy2}. 
\end{proof}

Compared with Theorem \ref{th:real}, Theorem \ref{th:general} is more restrictive in terms of the selection of the parameter $\gamma$. 
If $\Lambda_\tau H(z^\ast)$ has complex eigenvalues, then a smaller $\gamma$ is necessary to obtain a relation between asymptotically stable equilibria and saddle points. 
As outlined in the context of Theorem \ref{th:real}, for smaller $\eta$, the parameter $\gamma$ can be larger.
\begin{remark}
    In the absence of Assumption~\ref{as:inv}, 
    Theorems~\ref{th:real} and~\ref{th:general} hold
    only for strict saddle points.. 
    Without Assumption~\ref{as:inv}, we can not exclude the possibility that $\zeta = 0$ in proof of Theorem~\ref{th:real} and $(a,b)=(0,0)$ in proof of Theorem~\ref{th:general} at saddle points. 
    In these cases, the spectral radius may be equal to one, which does not necessarily imply asymptotic stability (contrary to the incorrect assertions 
    in 
    \cite{jin2020local,daskalakis2018limit}).
\end{remark}
Before we conclude 
this section, we show the application of Theorem \ref{th:general} based on an example also studied in the context of GDA algorithms in \cite{daskalakis2018limit}.
\begin{example} \label{ex:f_xy}
Consider the function $f:\R^2 \rightarrow\R$ defined as
\begin{align}
    f(x,y) =xy \label{eq:f_xy}
\end{align}
with unique critical point and unique saddle point $z^\ast=(x^\ast,y^\ast)=(0,0)$.
For $h_{1y}=\eta$ the Jacobian matrix \eqref{JGEG} is  	
\begin{align}
	&J(z^\ast) = 
	\begin{bmatrix}
		1-\frac{\gamma\eta^2}{\tau} & -\frac{\gamma\eta}{\tau}\\
		\gamma\eta &1-\frac{\gamma\eta^2}{\tau}
	\end{bmatrix}. \label{eq:J_xy}
\end{align}
with eigenvalues
 $(1-\frac{\gamma\eta^2}{\tau}) \pm
 \frac{\gamma\eta}{\sqrt{\tau}}\jmath$. 
 in Theorem~\ref{th:general}.  
For the parameter selection $\eta=0.9,$ $\tau=1,$ and $\gamma=0.1$ 
(satisfying the assumptions of Theorem \ref{th:general})
we have $\rho(J)<0.924$ and thus asymptotic stability of the equilibrium can be concluded from Theorem \ref{th:general}.

In \cite{daskalakis2018limit,jin2020local} stability properties of 
the dynamics of the $\tau$-GDA algorithm  are studied. The dynamics of $\tau$-GDA is described by  
\begin{align}
w_{\tau}(z) = z -\eta\Lambda_\tau F(z), \label{eq:GDA}
\end{align}
where $\eta>0$ is the step-size and $F$ and $\Lambda_\tau$ are defined in~\eqref{eq:F_definition}. 
The Jacobian of $\tau$-GDA  at critical points is $J_{\tau}(z)= I+\eta\Lambda_\tau H(z)$ for $H$  in \eqref{eq:DF_def}. 

Here, the eigenvalues of $H$ are purely complex, i.e., $\sigma(H(z^\ast)) =\{\pm\frac{1}{\tau} \jmath\}$, thus the results in \cite{jin2020local} (and in particular in \cite[Prop. 25]{jin2020local}) are not applicable to conclude that $(0,0)$ is a saddle point. 
Moreover, the eigenvalues of the Jacobian satisfy $\sigma(J_\tau(z^\ast))= \{1\pm\frac{\eta}{\tau} \jmath\}$, and thus $\rho(J_\tau(z^\ast))>1$ for all parameter selections of $\eta$ and $\tau$.
Hence, $z^\ast=(0,0)$ is unstable for \eqref{eq:GDA}. \hfill $\blacksquare$
\end{example}

Theorems \ref{th:real} and \ref{th:general} provide conditions such that the set of saddle points of a function $f$ are a subset of 
the set of asymptotically stable equilibria of \eqref{eq:discrete_time_sys}.
A natural follow up question is if there is a one-to-one mapping between saddle points and asymptotically stable equilibria. As the following counter example shows, the answer to this question generically is negative. 
Let $f(x,y) = -0.1x^2-0.5y^2+0.5xy$, which is a $C^3$-function that satisfies $\|\nabla^2f(x,y)\|\leq1.25$ for all 
$(x,y)\in \R^2$, and thus $\nabla f$ is globally Lipschitz with $L=1.25$ 
(satisfying the hypothesis of Theorem~\ref{th:general}).
The point 
$(0,0)$ is a critical point of $f$ which is not a saddle point  
due to the fact that $\nabla_{xx}^2 f<0$. 
Letting $h_{1y}=\eta=0.7$,
$\gamma=0.2$ and $\tau=1$ 
(satisfying the conditions of Theorem \ref{th:general}) it can be seen that $\rho(J(0,0))<1$, 
i.e.,
$(0,0)$ is asymptotically stable.

\section{Numerical example}\label{sec:exps}

\begin{table}[t]
\scriptsize
\centering
\setlength{\tabcolsep}{3pt}
\begin{tabular}{c|c|c|>{\centering\arraybackslash}p{1cm}|>{\centering\arraybackslash}p{1cm}|>{\centering\arraybackslash}p{1cm}}

\toprule
\multicolumn{3}{c|}{Parameters} & \multicolumn{3}{c}{Convergence (Asymptotic stability)} \\
\cmidrule(lr){1-3} \cmidrule(lr){4-6}
\(\eta\) & \(\tau\) & \(\gamma\) & $\tau$-EG & EG+ & GEG \\
\midrule
0.9 & 2 & 0.25 & YES & YES & YES \\
0.9 & 0.5 & 0.25 & NO & YES & YES \\
0.5 & 0.1 & 0.1 & NO & YES & YES \\
\midrule
0.9 & 2 & 1.2 & YES & NO & YES \\
0.9 & 0.01 & 0.1 & NO & YES & NO \\
0.9 & 0.01 & 0.01 & NO & YES & YES \\
0.5 & 2 & 2 & YES & NO & NO \\
0.5 & 200 & 2 & YES & NO & NO \\
\bottomrule
\end{tabular}
\caption{Convergence analysis of 3 variants of extra-gradient algorithm applied to Problem~\ref{eq:eq20} for \(f(x,y)=xy\).}\label{tab:xy}

\vspace{-0.4cm}
\end{table}%

We illustrate the results in Section \ref{sec:main} based on three examples.
First, we focus on $f(x,y)=xy$, a well-known example where GDA fails to converge (as discussed in Example \ref{ex:f_xy}).
Second, we construct a function with multiple critical points to analyse the convergence properties of \eqref{dy2} depending on the initial condition.
In the third example, a classifier neural network is trained and a corresponding empirical risk minimisation problem is solved using~\eqref{dy2}.

\subsection{Fixed points of GEG applied to $f(x,y) = xy$}

As the first example, we continue with the function $f$ discussed in Example \ref{ex:f_xy} 
and study the behaviours of three variants of extra-gradient algorithms: the $\tau$-EG algorithm (i.e., \eqref{dy2} for $\gamma=1$) \cite{chae2023two},  EG+ (i.e., \eqref{dy2} for $\tau=1$)  \cite{diakonikolas2021efficient}, and  GEG  introduced in this paper in \eqref{dy2}. 
We analyse the convergence of these three algorithms under different parameter selections. 
In Table~\ref{tab:xy}, in the first column the parameters of 
\eqref{dy}
are given. 
The $\tau$-EG algorithm has the same values for $\eta$ and $\tau$ as 
\eqref{dy}
(but $\gamma=1$) and the EG+ algorithm has the same values for $\eta$ and $\gamma$ as 
\eqref{dy}
(but $\tau=1$). 
The first three parameter selections in Table \ref{tab:xy} satisfy the assumptions of Theorem \ref{th:general}, and thus \eqref{dy} converges as expected. For the remaining five rows the Assumptions of Theorem \ref{th:general} are not satisfied, highlighting that our result is only sufficient, but not necessary.
Note that the dynamics for all the above three algorithms lead to a linear discrete-time system $z^+=Jz$ where $J$ defined in \eqref{eq:J_xy} is constant. 
Thus, in Table~\ref{tab:xy}, stability is verified through the eigenvalues of $J$.

\subsection{Fixed points of GEG applied to a function with multiple 
critical points}\label{sec:multiple}

\begin{figure}[b]
	\centering \includegraphics[width=0.45\textwidth]{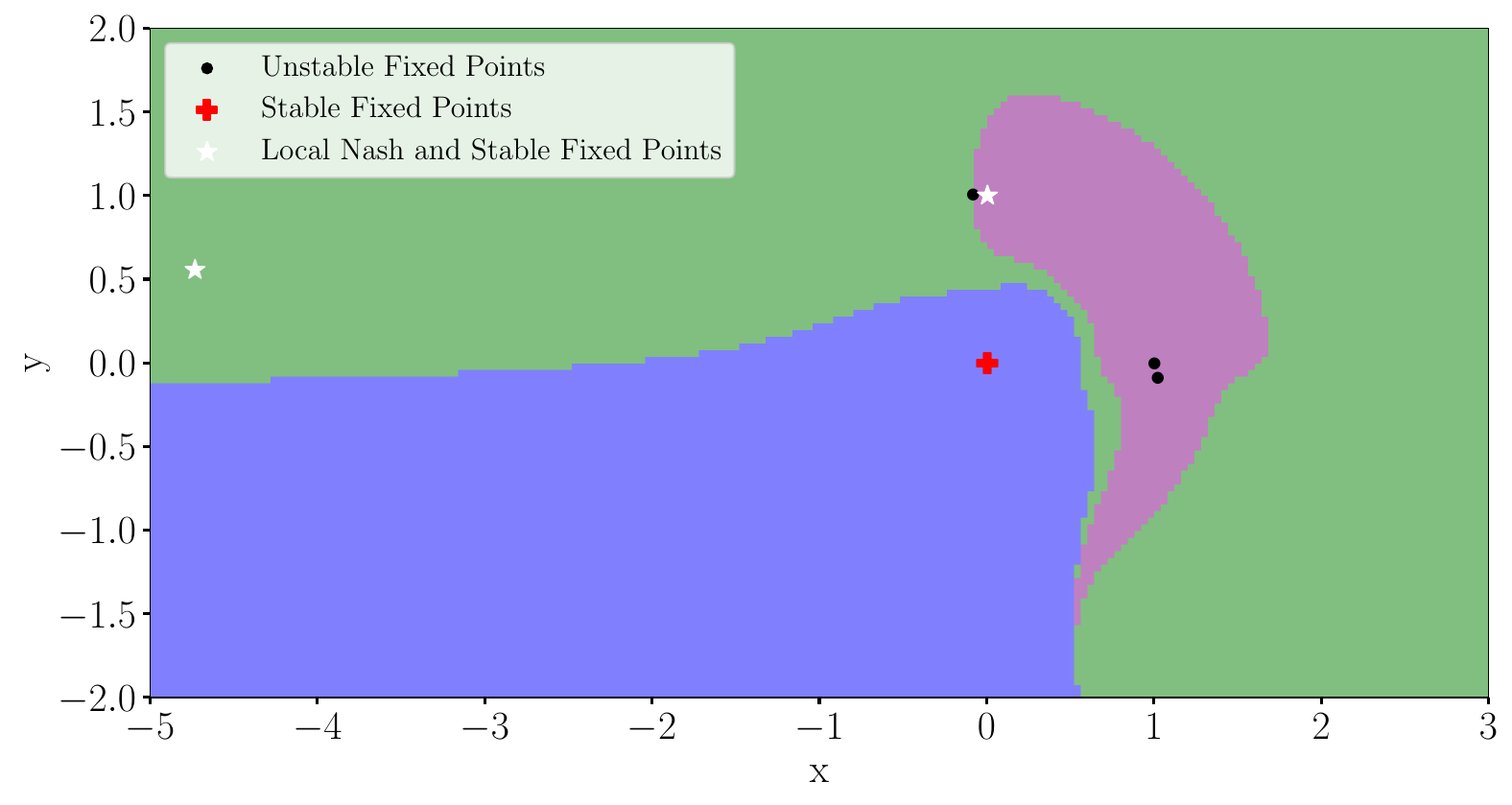}

\vspace{-0.4cm}
    
	\caption{Numerical estimates of regions
    of attraction of equilibria of the dynamics \eqref{dy2} for $f$ in 
    \eqref{eq:multiple}.}
	\label{fig:grid1}
\end{figure}

Following \cite[Sec 4.1]{daskalakis2018limit}, we construct a 
function with local saddle points, and with asymptotically stable and unstable equilibria with respect to \eqref{dy2}. 
Consider
\begin{align}\label{eq:multiple}
    f(x,y) &= f_1(x,y)(x-1)^2(y-1)^2+f_2(x,y)x^2y^2, 
\end{align}
where $f_1(x,y) = -0.25x^2 - 0.5y^2 + 0.6xy $ and 
    $f_2(x,y) = 0.5x^2 + 0.5y^2 + 4xy$, 
and with eight critical points reported in
Table~\ref{tab:conv}. The third column of Table \ref{tab:conv} states if a critical point is a saddle point of $f$ or not.
For the parameter selection $\eta\leq10^{-6},$ $\tau=1,$ $\gamma=0.5$,
satisfying the assumptions of Theorem \ref{th:general}, the second column indicates if a fixed point of the GEG-algorithm is stable. 
As expected from Theorem \ref{th:general}, the set of saddle points are a subset of the asymptotically stable points of \eqref{dy2}.
Additionally, estimates of the regions of attraction 
of three asymptotically stable equilibria on the domain $(x,y)\in [-5,3]\times[-2,2]$ for $\eta\leq10^{-4},$ $\tau=1,$ $\gamma=0.5$ are shown in Figure~\ref{fig:grid1}.
Here, the gradient of $f$ is locally Lipschitz but not globally Lipschitz, and thus, GEG renders the critical point $(38.402, -1.487)$ unstable if the step-size $\eta$ is not sufficiently small.
\begin{table}[htb]
\scriptsize	\centering
	\setlength{\tabcolsep}{4pt}
	\begin{tabular}{c|c|c}
		\toprule
        {Critical Points (Equilibria) } & GEG-stable & 
        Saddle point \\
		\midrule
		$(0,0)$& YES  & NO \\
		$(0,1)$& YES  & YES  \\
		$(1,0)$& NO  & NO \\
		$(-4.734, 0.560)$& YES  & YES   \\
		$(1.017, -0.086)$& NO  & NO \\
            $(0.731, -5.399)$ & NO & NO  \\
            $(-0.085, 1.006)$ & NO & NO \\
            $(38.402, -1.487)$ & YES & YES \\
		\bottomrule
	\end{tabular}
	\caption{Summary of critical points of \eqref{eq:multiple}.}
	\label{tab:conv}
\end{table}

\subsection{Large Scale Robust Optimisation} \label{sec:RO-Class}

\begin{figure}
	\centering \includegraphics[width=0.48\textwidth]{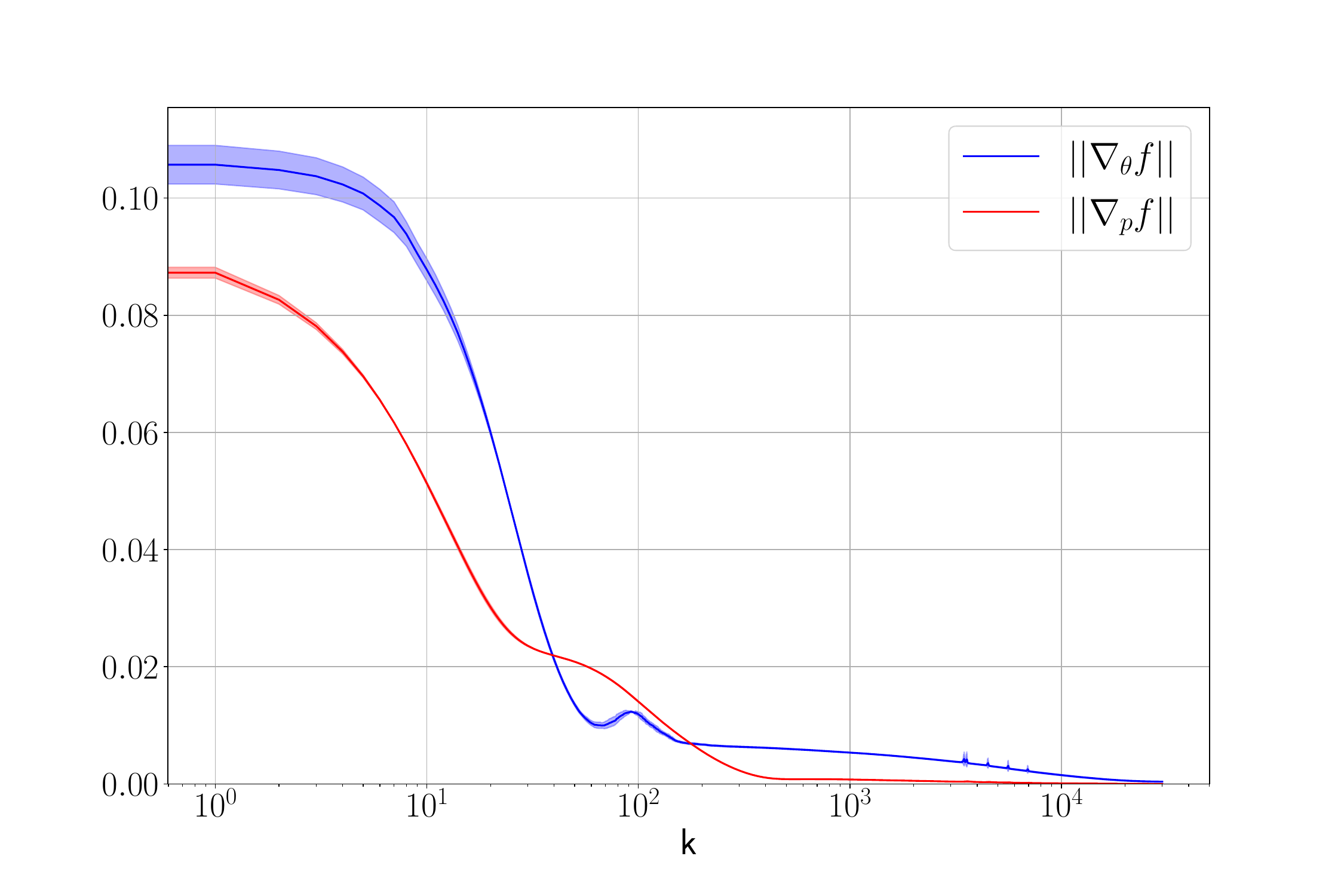}

    \vspace{-0.4cm}
    
    \caption{Evolution of gradient norm over \eqref{dy}} iterations.
	\label{fig:grad}
\end{figure}

For the third example, we consider a large scale optimisation problem whose critical points contrary to the previous two examples cannot be calculated explicitly. Specifically, we consider an
empirical risk minimisation for
a binary classification task given by 
\cite{anagnostidis2021direct}
\begin{align*} 
    &\underset{\theta}{\min}\;\underset{p}{\max}\; \textstyle
	-\sum_{i=1}^{m}p_i\big[y_i\log(\hat{y}(X_i;\theta))\\
    &\quad \textstyle +(1-y_i)\log(1-\hat{y}(X_i;\theta))\big]-\alpha\sum_{i=1}^{m}\left(p_i-\frac{1}{m}\right)^2.
\end{align*} 
Here,  $X_i\in \mathbb{R}^v$, $i\in \{1,\ldots,m\}$, are the data points, $\theta\in\mathbb{R}^n$ are the network parameters, and $\hat{y}(X_i;\theta),y_i \in \mathbb{R}^m$ are the predicted and the
true class of data  points 
$X_i$, respectively, and $p\in \mathbb{R}^m$ denotes the weights assigned to each data point.
The positive scalar $\alpha$ is the
regularisation parameter. 
We assume that $\hat{y}(X;\theta)$, is generated by a neural network with a hidden layer of size 50 and the LeakyReLU activation function with $n=1601$, and $m=455$ and where we use $80\%$ of 
the Wisconsin breast cancer data
set\footnote{\url{https://archive.ics.uci.edu/ml/datasets/Breast+Cancer+Wisconsin+(Diagnostic)}}
for training.
 Figure~\ref{fig:grad} shows the evolution of the norm of the gradient of a 5-fold cross-validation process \cite{Refaeilzadeh2009} of the GEG-algorithm using the parameters $\alpha=1,$ $\eta=0.01,$ $\tau=2,$ $\gamma=0.8,$ and randomly sampled initial values from a normal distribution.
The GEG-algorithm converges to a critical point of the objective function, and according to Theorem \ref{nouns} to asymptotically stable equilibria of the dynamics \eqref{dy2}.

\addtolength{\textheight}{-3cm}   

\section{Conclusions and Future Works}\label{sec:conc}
In this work, we introduced the Generalized Extra Gradient algorithm \eqref{dy}  and analyzed its limit points. We demonstrated that \eqref{dy2} avoids unstable critical points for almost all initialisations. From the dynamical systems perspective, we further proved that, under mild conditions, saddle points form a subset of asymptotically stable equilibria 
of \eqref{dy2}.
For future work, we plan to explore the connection between local min-max points, as defined in \cite{jin2020local}, and GEG-stable points. Additionally, we aim to investigate the relation between the saddle points and the continuous-time dynamics. 

\bibliographystyle{ieeetr}  

\bibliography{EG_stability1.bib}

\begin{thebibliography}{10}

\bibitem{myerson2013game}
R.~B. Myerson, {\em Game Theory}.
\newblock Harvard University Press, 2013.

\bibitem{madry2018towards}
A.~Madry, A.~Makelov, L.~Schmidt, D.~Tsipras, and A.~Vladu, ``Towards deep
  learning models resistant to adversarial attacks,'' in {\em International
  Conference on Learning Representations}, 2018.

\bibitem{omidshafiei2017deep}
S.~Omidshafiei, J.~Pazis, C.~Amato, J.~P. How, and J.~Vian, ``Deep
  decentralized multi-task multi-agent reinforcement learning under partial
  observability,'' in {\em International Conference on Machine Learning},
  pp.~2681--2690, PMLR, 2017.

\bibitem{korpelevich1977extragradient}
G.~M. Korpelevich, ``Extragradient method for finding saddle points and other
  problems,'' {\em Matekon}, vol.~13, no.~4, pp.~35--49, 1977.

\bibitem{chen1997convergence}
G.~H. Chen and R.~T. Rockafellar, ``Convergence rates in forward--backward
  splitting,'' {\em SIAM Journal on Optimization}, vol.~7, no.~2, pp.~421--444,
  1997.

\bibitem{wei2021last}
C.-Y. Wei, C.-W. Lee, M.~Zhang, and H.~Luo, ``Last-iterate convergence of
  decentralized optimistic gradient descent/ascent in infinite-horizon
  competitive markov games,'' in {\em Conference on learning theory},
  pp.~4259--4299, PMLR, 2021.

\bibitem{beznosikov2023stochastic}
A.~Beznosikov, E.~Gorbunov, H.~Berard, and N.~Loizou, ``Stochastic gradient
  descent-ascent: Unified theory and new efficient methods,'' in {\em
  International conference on artificial intelligence and statistics},
  pp.~172--235, PMLR, 2023.

\bibitem{daskalakis2018limit}
C.~Daskalakis and I.~Panageas, ``The limit points of (optimistic) gradient
  descent in min-max optimization,'' {\em Advances in neural information
  processing systems}, vol.~31, 2018.

\bibitem{goebel2024set}
R.~K. Goebel, {\em Set-Valued, Convex, and Nonsmooth Analysis in Dynamics and
  Control: An Introduction}.
\newblock SIAM, 2024.

\bibitem{chae2023two}
J.~Chae, K.~Kim, and D.~Kim, ``Two-timescale extragradient for finding local
  minimax points,'' in {\em The Twelfth International Conference on Learning
  Representations, ICLR 2024}, International Conference on Learning
  Representations, 2024.

\bibitem{diakonikolas2021efficient}
J.~Diakonikolas, C.~Daskalakis, and M.~I. Jordan, ``Efficient methods for
  structured nonconvex-nonconcave min-max optimization,'' in {\em International
  Conference on Artificial Intelligence and Statistics}, pp.~2746--2754, PMLR,
  2021.

\bibitem{kim2024double}
K.~Kim and D.~Kim, ``Double-step alternating extragradient with increasing
  timescale separation for finding local minimax points: provable
  improvements,'' in {\em Proceedings of the 41st International Conference on
  Machine Learning}, pp.~24059--24093, 2024.

\bibitem{hespanha2018linear}
J.~P. Hespanha, {\em Linear Systems Theory}.
\newblock Princeton University Press, 2018.

\bibitem{kellett2023introduction}
C.~M. Kellett and P.~Braun, {\em Introduction to Nonlinear Control: Stability,
  Control Design, and Estimation}.
\newblock Princeton University Press, 2023.

\bibitem{bertsekas2009convex}
D.~Bertsekas, {\em Convex Optimization Theory}.
\newblock Athena Scientific, 2009.

\bibitem{jin2020local}
C.~Jin, P.~Netrapalli, and M.~Jordan, ``What is local optimality in
  nonconvex-nonconcave minimax optimization?,'' in {\em International
  conference on machine learning}, pp.~4880--4889, PMLR, 2020.

\bibitem{lee2012smooth}
J.~M. Lee, {\em Smooth Manifolds}, pp.~1--31.
\newblock New York, NY: Springer New York, 2012.

\bibitem{rudin1964principles}
W.~Rudin {\em et~al.}, {\em Principles of Mathematical Analysis}, vol.~3.
\newblock McGraw-hill New York, 1964.

\bibitem{moslehian2012ky}
M.~S. Moslehian, ``Ky fan inequalities,'' {\em Linear and multilinear algebra},
  vol.~60, no.~11-12, pp.~1313--1325, 2012.

\bibitem{lee2019first}
J.~D. Lee, I.~Panageas, G.~Piliouras, M.~Simchowitz, M.~I. Jordan, and
  B.~Recht, ``First-order methods almost always avoid strict saddle points,''
  {\em Mathematical programming}, vol.~176, pp.~311--337, 2019.

\bibitem{anagnostidis2021direct}
S.-K. Anagnostidis, A.~Lucchi, and Y.~Diouane, ``Direct-search for a class of
  stochastic min-max problems,'' in {\em International Conference on Artificial
  Intelligence and Statistics}, pp.~3772--3780, PMLR, 2021.

\bibitem{Refaeilzadeh2009}
P.~Refaeilzadeh, L.~Tang, and H.~Liu, {\em Cross-Validation}, pp.~532--538.
\newblock Boston, MA: Springer US, 2009.

\end{thebibliography}


\appendix

\section{Mathematica Code used in the proof of Theorem~\ref{th:general}} 
\subsection{ Mathematica Code for bounding $\gamma$ such that $|\lambda|\leq1$}\label{mtc}
\begin{lstlisting}	
f[a_, b_, \[Gamma]_] := 
 Sqrt[(1 + \[Gamma] (a + a^2 - b^2))^2 + \[Gamma]^2 (b + 2 a b)^2]
constraints = -1 < a < 0 && -Sqrt[1 - a^2] < b < Sqrt[1 - a^2];
maxConstraint = MaxValue[{f[a, b, \[Gamma]], constraints}, {a, b}];
\[Gamma]Bounds = Reduce[maxConstraint <= 1, \[Gamma]]
\end{lstlisting}
\subsection{Mathematica Code for finding a feasible nonzero $\gamma$}\label{mtc2}
\begin{lstlisting}	
f[a_, b_] := -2 (a^2 - b^2 + a)/((a^2 - b^2 + a)^2 + (2 a b + b)^2)
domain = -1 < a <= 0 && -Sqrt[1 - a^2] < b < Sqrt[1 - a^2] && {a, b} != {0, 0};
MinValue[{f[a, b], domain}, {a, b}]
ArgMin[{f[a, b], domain}, {a, b}]
MaxValue[{f[a, b], domain}, {a, b}]
\end{lstlisting}

\end{document}